\documentclass[10pt, reqno]{amsart}
\usepackage{amsmath}
\usepackage{amsfonts,amssymb,amscd,amsmath, graphicx, bbm}

\newtheorem{theorem*}{Theorem}
\makeatletter
\def\cfK{%
    \operatornamewithlimits{%
        \mathchoice{
            \vcenter{\hbox{\huge $\mathbf{K}$}}%
        }{
            \vcenter{\hbox{\Large $\mathbf{K}$}}%
        }{
            \mathrm{K}%
        }{
            \mathrm{K}%
        }
    }
}
\makeatother
\begin{document}

\title{Some continued fractions for $\pi$ and $G$}
\date{\today}
\author {Amrik Singh Nimbran}
\address{B3-304, Palm Grove Heights, Ardee City, Gurgaon, Haryana, INDIA }
\email{amrikn622@gmail.com}

\author {Paul Levrie}
\address{Faculty of Applied Engineering, University of Antwerp, Groenenborgerlaan 171, B-2020 Antwerpen;
Department of Computer Science, KU Leuven, P.O. Box 2402, B-3001 Heverlee, BELGIUM}
\email{paul.levrie@cs.kuleuven.be}


\begin{abstract}
We present here two classes of infinite series and the associated continued fractions involving $\pi$ and  Catalan's constant based on the work of Euler and Ramanujan. A few sundry continued fractions are also given.
\end{abstract}

\maketitle

\section{Introduction}

A continued fraction is an expression of the general form

\[
b_0+\cfrac{a_1}{b_1+ \cfrac{a_2}{b_2+ \cfrac{a_3}{b_3+ \dotsb}}}
\]~\\
with $a_1,a_2,\ldots$ and $b_0, b_1, b_2,\ldots$ real numbers. The following space-saving notation can be used:
\[
\newcommand{\cfplus}{\mathbin{\genfrac{}{}{0pt}{}{}{+}}}
b_0+\frac{a_1}{b_1}\cfplus\frac{a_2}{b_2}\cfplus\frac{a_3}{b_3}\cfplus\dotsb
\]
or the shorter notation
$ b_0 +\cfK_{n=1}^\infty \displaystyle\frac{a_n}{b_n},$
where the letter K comes from the German word \emph{Kettenbr\"{u}che} (Ketten - chain, Br\"{u}che - fraction) for a continued fraction.

It is not known when, where and by whom continued fractions were first used, but the notion seems to be quite old. For an overview, see \cite{brezinski}. Euclid's algorithm for finding the greatest common devisor of two integers in effect converts a fraction into a terminating continued fraction. Continued fractions were used in India in the 6th century by Aryabhata to solve linear Diophantine equations and in the 12th century by Bhaskaracharya for solving the `Pell' equation. Rafael Bombelli (1526--1572) and Pietro Antonio Cataldi (1548--1626), both of Bologna (Italy), gave continued fractions for $\sqrt{13}$ and $\sqrt{18}$ respectively. The first continued fraction for a number other than the quadratic irrationals is due to William Brouncker (1620-84) who converted the product $\displaystyle \frac{4}{\pi}=\prod_{n=1}^\infty \frac{(2n +1)^2}{2n(2n+2)}$ submitted to him by John Wallis (1616--1703) in 1655 (see \cite[p.182]{wallis}):
\[
\frac{4}{\pi}=1 +\cfK_{n=1}^\infty \frac{(2n-1)^2}{2}=1+\cfrac{1^2}{2+ \cfrac{3^2}{2+ \cfrac{5^2}{2+ \dotsb}}} .
\]

It was one of the sequence of continued fractions found by Brouncker.\cite[p.307]{stedall}

\section{Euler's continued fractions for $\pi$}

Leonhard Euler (1707--1783) laid the foundation of the theory of continued fractions in a dozen papers (six published in his life-time and six posthumously) and chapter 18 of his \textit{Introductio in Analysin Infinitorum} (1748). Of these, five papers written between 1737 and 1780, with Enestr\"{o}m index numbers: E071, E123, E522, E593 and E745, are of interest to us here. Brouncker's fraction occurs in E071 \cite[\S 4]{euler71}.
Euler explains the method of converting an infinite series into a continued fraction in E593 \cite{euler593}. He gives among others this formula in \cite[\S7]{euler123}:

\[
\int_{0}^1 \frac{x^{n-1}{dx}}{1+x^m}=\cfrac{1}{n+ \cfrac{n^2}{m+ \cfrac{(m+n)^2}{m +\cfrac{(2m+n)^2}{m+ \cfrac{(3m+n)^2}{m+ \dotsb}}}}} .
\]

Now $\displaystyle \int \frac{x {dx}}{1+x^4}=\frac{1}{2}\arctan x^2 +C.$ So taking $m=4, \, n=2,$ we get Brouncker's fraction in disguise:

\[
\frac{\pi}{8}=\cfrac{1}{2+ \cfrac{2^2}{4+ \cfrac{6^2}{4 +\cfrac{10^2}{4+ \cfrac{14^2}{4+ \dotsb}}}}} .
\]

We find in \cite[\S18]{euler593} Euler's Theorem I:

\textit{If such an infinite series will have been proposed
\[
s=\frac{1}{\alpha} -\frac{1}{\beta} +\frac{1}{\gamma} -\frac{1}{\delta} +\frac{1}{\epsilon} -\cdots,
\]
a continued fraction of the form
\[
\frac{1}{s}=\alpha +{\cfrac{\alpha \alpha}{\beta-\alpha+ \cfrac{\beta\beta}{\gamma-\beta+ \cfrac{\gamma\gamma}{\delta -\gamma+ \dotsb}}}}
\]
can always be formed from it.}

Taking the Leibnitz series $\displaystyle \sum_{n=1}^\infty \frac{(-1)^{n-1}}{2n-1}=\frac{\pi}{4},$ Euler deduced Brouncker's fraction.

We find this theorem as Theorem II in \cite[\S23]{euler593}:

\textit{If the proposed series is of this form
\[
s=\frac{1}{ab} -\frac{1}{bc} +\frac{1}{cd} -\frac{1}{de} +\frac{1}{ef} -\cdots,
\]
from this, the following continued fraction springs forth
\[
\frac{1}{a s}=b +{\cfrac{ab}{c-a+ \cfrac{bc}{d-b+ \cfrac{cd}{e -c+\cfrac{de}{f-d \dotsb}}}}} .
\]
}
Noting that $\displaystyle \sum_{n=1}^\infty \frac{(-1)^{n-1}}{(2n-1)(2n+1)}=\frac{\pi}{4} -\frac{1}{2},$ he deduced
\[
\frac{4}{\pi -2}=3 +{\cfrac{1\cdot3}{4+ \cfrac{3\cdot5}{4+ \cfrac{5\cdot7}{4+ \cfrac{7\cdot9}{4+\dotsb}}}}}
\]
or
\[
\frac{\pi}{2}-1 =\cfrac{2}{3+ \cfrac{1\cdot3}{4+ \cfrac{3\cdot5}{4+ \cfrac{5\cdot7}{4+ \dotsb}}}} .
\]

Note that the convergents of this continued fraction $\displaystyle c_n=\frac{p_n}{q_n}$ satisfy the following relations:
\[ p_0=q_0=1, \quad
p_{n+1}=(2n+3) p_n +(-1)^n 2 (2n-1)!!, \quad q_n=(2n+1)!!,
\]
for $n=0, 1, 2, \ldots$,
where $n!!=n \cdot (n-2) \cdot (n-4) \cdot \ldots \cdot 1$ for $n$ odd, with $(-1)!!=1$.

This continued fraction is a special case of a more general one which can be deduced from a formula of Ramanujan giving a continued fraction for a product of quotients of gamma function values \cite[p.227]{perron}:

\noindent{\textbf{General Formula}}: \textit{Let
\[
P_{-1}=\frac{1}{2}, \, \ P_{2m}=\prod_{k=1}^{m} \frac{2k(2k+2)}{(2k+1)^2}, \, \  P_{2m+1}=P_{2m}\cdot \frac{2m+2}{2m+3} \quad (m=0,1,2,\ldots).
\]
Then we have
\[
\frac{\pi}{4 P_n} -1 =\cfrac{(-1)^{n+1} 2}{4(n+2)+(-1)^n+ \cfrac{1\cdot3}{4(n+2)+ \cfrac{3\cdot5}{4(n+2)+ \cfrac{5\cdot7}{4(n+2) +\dotsb}}}}
\]
}

It may be compared with \cite[p.227, eq(30a)]{perron}. Note that these continued fractions give $\pi$ as a combination of the partial products of Wallis's formula and a continued fraction. For instance, for $n=4$ we have:
\[
\frac{\pi}{4}=\frac{2}{3}\cdot\frac{4}{3}\cdot\frac{4}{5}\cdot\frac{6}{5} \cdot   [ 1- \cfrac{2}{25+ \cfrac{1\cdot3}{24+ \cfrac{3\cdot5}{24+ \cfrac{5\cdot7}{24+\ddots}}}} ] .
\]

Euler gave the following two continued fractions and a couple more in \cite[\S31, \S33]{euler123}:
\[
\frac{\pi}{2}=1 +\cfrac{1}{1+ \cfrac{1\cdot2}{1+ \cfrac{2\cdot3}{1+ \cfrac{3\cdot4}{1+ \dotsb}}}} \ ,
\]
\[
\frac{\pi}{2}=2- \cfrac{1}{2+ \cfrac{1^2}{2+ \cfrac{2^2}{2+ \cfrac{3^2}{2+ \dotsb}}}} .
\]

He obtained this fraction with partial denominators $4n$ from his general formula in \cite[\S 18]{euler522} and \cite[\S36]{euler593}:
\[
\frac{6\sqrt{3}} {\pi} =3 +\cfrac{3\cdot1^2}{8+\cfrac{3\cdot3^2}{12+ \cfrac{3\cdot5^2}{16+ \cfrac{3\cdot7^2} {20+ \dotsb}}}} .
\]

These continued fractions with partial denominators $3n-2$ and $5n-3$ are due to Glaisher \cite{glaisher73, glaisher76}.

\[
\frac{2}{\pi}=1-\cfrac{1\cdot1}{4- \cfrac{2\cdot3}{7- \cfrac{3\cdot5}{10- \cfrac{4\cdot7} {13- \dotsb}}}} \ ,
\]
\[
\frac{3\sqrt{3}}{\pi}=2-\cfrac{2(1\cdot1)}{7- \cfrac{2(2\cdot3)}{12- \cfrac{2(3\cdot5)}{17- \dotsb}}} .
\]
They were converted from the following series:
\[
\frac{\pi}{2}=1 +\frac{1}{3} +\frac{1\cdot2}{3\cdot5} +\frac{1\cdot2\cdot3}{3\cdot5\cdot7} +\frac{1\cdot2\cdot3\cdot4}{3\cdot5\cdot7\cdot9}+\cdots\ ,
\]
\[
\frac{2\pi}{3\sqrt{3}}=1 +\frac{1}{2\cdot3} +\frac{1\cdot2}{3\cdot4\cdot5} +\frac{1\cdot2\cdot3}{4\cdot5\cdot6\cdot7} +\frac{1\cdot2\cdot3\cdot4}{5\cdot6\cdot7\cdot8\cdot9}+\cdots.
\]

We will now obtain some new continued fractions for $\pi.$

\section{Series with linear factors and continued fractions}

If we define:
\[
y_k = \sum_{n=1}^\infty \frac{(-1)^{n-1}}{(2n-1)(2n+1)(2n+3)\cdots(2n+2k-3)(2n+2k-1)}
\]
then the following recurrence relation with $k \in \mathbb{N}$ is easy to establish:
\[
y_k = \frac{1}{k} \cdot y_{k-1} -\frac{1}{2k \cdot 1\cdot3\cdot5\cdots(2k-3)\cdot (2k-1)}.
\]
Indeed, we have that:
\begin{align*}
2k y_ k & = 2k \sum_{n=1}^\infty \frac{(-1)^{n-1}}{(2n-1)(2n+1)\cdots(2n+2k-3)(2n+2k-1)}\\
& = \sum_{n=1}^\infty \frac{(-1)^{n-1}}{(2n-1)(2n+1)\cdots(2n+2k-5)(2n+2k-3)} \\
& \hspace{4cm} -\sum_{n=1}^\infty \frac{(-1)^{n-1}}{(2n+1)\cdots(2n+2k-3)(2n+2k-1)} .
\end{align*}
By shifting the index of summation, the last sum on the right can be rewritten as:
\begin{align*}
& -\sum_{n=2}^\infty \frac{(-1)^{n-1}}{(2n-1)(2n+1)\cdots(2n+2k-5)(2n+2k-3)}\\
& = \frac{1}{1\cdot3\cdot5\cdots(2k-1)} -\sum_{n=1}^\infty \frac{(-1)^{n-1}}{(2n-1)(2n+1)\cdots(2n+2k-5)(2n+2k-3)}
\end{align*}
immediately leading to the expected result. See also \cite{nimbrie}.

Using this recurrence with $y_0 = \frac{\pi}{4}$ (Leibnitz's series), we get a class of series with an increasing number of factors in the denominator:
\begin{equation}
\frac{\pi}{8}-\frac{1}{3}=\sum_{n=1}^\infty \frac{(-1)^{n-1}}{(2n-1)(2n+1)(2n+3)} \ , \label{eq:a}
\end{equation}

\[
\frac{\pi}{24}-\frac{11}{90}=\sum_{n=1}^\infty \frac{(-1)^{n-1}}{(2n-1)(2n+1)(2n+3)(2n+5)} \ ,
\]

\[
\frac{\pi}{96}-\frac{2}{63}=\sum_{n=1}^\infty \frac{(-1)^{n-1}}{(2n-1)(2n+1)(2n+3)(2n+5)(2n+7)} .
\]

Following Euler's method, we obtain
\[
\frac{8}{3\pi -8}=5 +{\cfrac{1\cdot5}{6+ \cfrac{3\cdot7}{6+ \cfrac{5\cdot9}{6+ \cfrac{7\cdot11}{6+ \dotsb}}}}} \ ,
\]

\[
\frac{24}{15\pi -44}=7 +{\cfrac{1\cdot7}{8+ \cfrac{3\cdot9}{8+ \cfrac{5\cdot11}{8+ \cfrac{7\cdot13}{8+ \dotsb}}}}} \ ,
\]

\[
\frac{96}{105\pi -320}=9 +{\cfrac{1\cdot9}{10+ \cfrac{3\cdot11}{10+ \cfrac{5\cdot13}{10+ \cfrac{7\cdot15}{10+ \dotsb}}}}} .
\]

It may be noted that the numerator of the left hand side is $4(k!)$ while the multiple of $\pi$ is $(2k-1)!!,$ where $(2k-1)$ is the constant in the largest factor $2n+2k-1$ in the denominator of the series. These continued fractions are all special cases of:

\begin{theorem*}
Let $s=\displaystyle \sum_{n=1}^\infty \frac{(-1)^{n-1}}{\prod_{j=1}^k{(2n+2j-3)}},$ then we have
\[
\frac{1}{(2k-1)!! \, s}=2k+1 +{\cfrac{1\cdot(2k+1)}{2k+2+ \cfrac{3\cdot(2k+3)}{2k+2+ \cfrac{5\cdot(2k+5)}{2k+2+ \cfrac{7\cdot(2k+7)}{2k+2+ \dotsb}}}}} \ .
\]
\end{theorem*}
To illustrate how the continued fractions are derived from the corresponding series (Euler's theorem I), we give an example (see also \cite{levrie}).
We will convert the following series given in \cite[p.269, Ex.109(c)]{knopp}:
\begin{equation}
\sum_{n=1}^\infty \frac{(-1)^{n-1}}{2n(2n+1)(2n+2)}=\frac{\pi-3}{4} \label{eq:b}
\end{equation}
into
\[
\frac{1}{\pi -3}=6 +{\cfrac{1^2}{6+ \cfrac{3^2}{6+ \cfrac{5^2}{6+ \dotsb}}}} ,
\]
the ``new" continued fraction derived in 1999 by Lange \cite{lange} though it already occurs in a 1988 paper of Castellanos \cite{castellanos} who rightly ascribed it to Euler. This continued fraction is related to Brouncker's.

Let us define $z_n = (-1)^n (\frac{\pi-3}{4}-\sum_{k=1}^n \frac{(-1)^{k-1}}{2k(2k+1)(2k+2)})$. From this it follows that:
\[
z_{n}+z_{n-1} =\frac{1}{2n(2n+1)(2n+2)} , \ \mbox{and} \
z_{n+1}+z_{n} =\frac{1}{(2n+2)(2n+3)(2n+4)} .
\]
Dividing both equations, we find:
\[
\frac{z_{n+1}+z_n}{z_n+z_{n-1}}=\frac{n(n+1)(2n+1)}{(n+1)(n+2)(2n+3)}=\frac{n(2n+1)}{(n+2)(2n+3)} .
\]

We now get rid of the denominators:
\[
(2n+3)(n+2) z_{n+1}+[(2n+3)(n+2)-(2n+1) n]z_n = (2n+1) n z_{n-1}
\]
or
\[
(2n+3)(n+2) z_{n+1}+6 (n+1) z_n = (2n+1) n z_{n-1} .
\]

If we use the transformation $w_n = (n+1) z_n$, we have:
\[
(2n+3) w_{n+1}+6w_n = (2n+1) w_{n-1} \Rightarrow (2n+3)\frac{w_{n+1}}{w_n} +6=\frac{2n+1}{\displaystyle \frac{{w_n}}{w_{n-1}}}
\]
giving the relation that generates the continued fraction on taking $n=1, 2, 3, \cdots$
\[
\frac{{w_n}}{w_{n-1}}=\frac{(2n+1)}{6 +(2n+3)\displaystyle \frac{w_{n+1}}{w_n}}.
\]
The value of the continued fraction obtained in this way is given by:
\[
\frac{w_1}{w_0} = \frac{2z_1}{z_0}=2 \cdot \frac{- (\frac{\pi-3}{4}-\frac{1}{2\cdot 3\cdot 4})}{\frac{\pi-3}{4}} .
\]
Some manipulations lead to the desired form.
Other continued fractions for $\pi$ can be obtained in a similar way. We derived this interesting series by combining (\ref{eq:a}) and (\ref{eq:b}):

\begin{equation}
\sum_{n=1}^\infty \frac{(-1)^{n-1}}{(2n-1)2n(2n+1)(2n+2)(2n+3)}=\frac{10-3\pi}{72}
\end{equation}
which gives the second convergent for $\pi:$
\[
\pi=\frac{22}{7} -24\sum_{n=2}^\infty \frac{(-1)^{n}}{(2n+1)(2n+2)(2n+3)(2n+4)(2n+5)}.
\]

We converted the preceding series into the following continued fraction:

\[
\frac{6}{10 -3\pi}=10 +{\cfrac{1\cdot5}{10+ \cfrac{3\cdot7}{10+ \cfrac{5\cdot9}{10+ \cfrac{7\cdot11}{10+ \cfrac{9\cdot13}{10+ \dotsb}}}}}}
\]
which may be compared with the continued fraction derived earlier (with $a_k=10, \; \forall \; k\in \mathbb{N}$) and one given by Euler in \cite[\S12]{euler745} and also by Osler \cite{osler}:

\[
\frac{16}{\pi}=5 +{\cfrac{1^2}{10+ \cfrac{3^2}{10+ \cfrac{5^2}{10+ \cfrac{7^2}{10+ \cfrac{9^2}{10+ \dotsb}}}}}}
\]

Osler \cite{osler} (also in \cite[p.226, eq(28a)]{perron}) gives the following two classes of continued fractions related to Brouncker's. They were already known to John Wallis \cite{wallis}. Let
\[
P_0=1, \ \mbox{and} \ P_n=\prod_{k=1}^n \frac{(2k-1)(2k+1)}{(2k)^2}.
\]
Then for $n=0, 1, 2, 3, \cdots$ we have:
\begin{gather*}
(4n+1) +\cfrac{1^2}{2(4n+1)+ \cfrac{3^2}{2(4n+1)+ \cfrac{5^2}{2(4n+1)+ \dotsb}}}=\frac{(2n+1)}{P_n}\frac{4}{\pi},\\
(4n+3) +\cfrac{1^2}{2(4n+3)+ \cfrac{3^2}{2(4n+3)+ \cfrac{5^2}{2(4n+3)+ \dotsb}}}=(2n+1)P_{n} \pi.
\end{gather*}

These formulas are special cases of the known formula
\[
\frac{4 \, \displaystyle \Gamma\left(\frac{x+3 +y}{4}\right) \Gamma\left(\frac{x+3 -y}{4}\right)}{\displaystyle \Gamma\left(\frac{x+1 +y}{4}\right) \Gamma\left(\frac{x+1- y}{4}\right)}=x +\cfrac{1^2 -y^2}{2x+ \cfrac{3^2 -y^2}{2x+ \cfrac{5^2 -y^2}{2x+ \dotsb}}}
\]
valid for either $y$ an odd integer and $x$ any complex number or $y$ any complex number and $\Re(x)>0.$ Originally due to Euler \cite[\S67]{euler123}, it occurs in an inverted form in \cite[p.140]{bernt2}.

\section{Series with quadratic factors and continued fractions}

Ramanujan's Notebook II \cite[pp.151, 153]{bernt2} contains these two continued fractions for  Catalan's constant $\displaystyle G=\sum_{n=1}^\infty \frac{(-1)^{n-1}}{(2n-1)^2}=0.915965594177\dots$:

\[
2G=2 -{\cfrac{1}{3+ \cfrac{2^2}{1+ \cfrac{2^2}{3+ \cfrac{4^2}{1+ \cfrac{4^2}{3+ \dotsb}}}}}}
\]
and
\[
2G=1 +{\cfrac{1^2}{\frac{1}{2}+ \cfrac{1\cdot2}{\frac{1}{2}+ \cfrac{2^2}{\frac{1}{2}+ \cfrac{2\cdot3}{\frac{1}{2}+ \cfrac{3^2}{\frac{1}{2}+ \dotsb}}}}}}
\]

We find in Ramanujan's Manuscript Book 1 \cite[Ch.XIV, p.107, entry 14]{Ramanujan1} (and in \cite[p.123, entry 16]{bernt2}) this formula with $m, n \notin \mathbb{Z^{-}}$:

\[
\sum_{k=1}^\infty \frac{(-1)^{k+1}}{(m+k)(n+k)}=\cfrac{1}{m+n+1+mn+ \cfrac{(m+1)^2 (n+1)^2}{m+n+3+ \cfrac{(m+2)^2 (n+2)^2}{m+n+5+ \cfrac{(m+3)^2 (n+3)^2}{m+n+7+ \dotsb}}}}
\]

Setting $m=n=-\frac{1}{2}$ and doing a little manipulation, the last formula yields a fraction given in \cite{bowman}:
\begin{equation}
G=\cfrac{1}{1+ \cfrac{1^4}{8+ \cfrac{3^4}{16+ \cfrac{5^4} {24+ \dotsb}}}} .\label{eq:bowman}
\end{equation}
Using a generalisation of the series defining $G$, we will obtain new continued fractions for this constant and for $\pi$. In the \textit{Appendix} we prove that $y_k$ defined by
\[
y_k = \sum_{n=1}^\infty \frac{(-1)^{n-1}}{(2n-1)^2(2n+1)^2(2n+3)^2 \cdots (2n+2k-3)^2}
\]
satisfies the recurrence
\[
y_k=\frac{10k^2+8k+1}{2(2k+1)!!^2} -4k(k+1)^3 y_{k+2} .
\]
Since $y_1 = G$ and $y_2 = \frac{1}{2} - \frac{\pi}{8}$ (which can be proved using a partial fraction expansion and telescoping), the previous recurrence gives for $y_4$, $y_6$ and $y_8$:

\[
\sum_{n=1}^\infty \frac{(-1)^{n-1}}{(2n-1)^2(2n+1)^2(2n+3)^2(2n+5)^2}=\frac{2!^3}{4!^3 1!}\pi -\frac{7}{4050}
\]

\[
\sum_{n=1}^\infty \frac{(-1)^{n-1}}{(2n-1)^2(2n+1)^2\cdots(2n+9)^2}=-\frac{3!^3}{6!^3 2!}\pi +\frac{41}{44651250}
\]

\[
\sum_{n=1}^\infty \frac{(-1)^{n-1}}{(2n-1)^2(2n+1)^2\cdots(2n+13)^2}=\frac{4!^3}{8!^3 3!}\pi -\frac{14789}{134221791453750}
\]
and for $y_3$, $y_5$ and $y_7$:
\[
\sum_{n=1}^\infty \frac{(-1)^{n-1}}{(2n-1)^2(2n+1)^2(2n+3)^3}=-\frac{1!!^3 1!}{2!^4 2^1} G + \frac {19}{576}
\]
\[
\sum_{n=1}^\infty \frac{(-1)^{n-1}}{(2n-1)^2(2n+1)^2(2n+3)^2(2n+5)^2(2n+7)^3}=\frac{3!!^3 2!}{4!^4 2^2} G - \frac {3919}{108380160}
\]
\[
\sum_{n=1}^\infty \frac{(-1)^{n-1}}{(2n-1)^2(2n+1)^2\cdots(2n+11)^2}=-\frac{5!!^3 3!}{6!^4 2^3} G + \frac {22133579}{2549361475584000} .
\]
Note that an \emph{odd number} of consecutive odd squares in the denominator gives a series for $G$ while an \emph{even number} of these squares gives a series for $\pi.$

In general we found the following sums for the series related to $\pi$: $(k=0, 1, 2, \cdots)$
\[
\sum_{n=1-k}^\infty \frac{(-1)^{n-1}}{\prod_{j=0}^{2k+1} (2n+2j-1)^2}=(-1)^{k+1} \left(\frac{(k+1)!^3}{(2k+2)!^3 k!} \pi - \frac{1}{2(2k+1)!!^4}\right).
\]

For the series above the corresponding continued fractions are given in the following theorem:
\begin{theorem*}
Let
\[S
=\displaystyle \sum_{n=1}^\infty \frac{(-1)^{n-1}}{(2n-1)^2(2n+1)^2(2n+3)^2 \cdots (2n+(2k-1))^2}.\]
Then,
\[
 \frac{1}{(2k-1)!!^2 \, S}=(2k+1)^2 +\cfrac{1^2 (2k+1)^2}{(2k+3)^2-1^2+ \cfrac{3^2 (2k+3)^2}{(2k+5)^2-3^2+ \cfrac{5^2 (2k+5)^2}{(2k+7)^2-5^2+ \dotsb}}} .
\]
\end{theorem*}

Using the above formulas and the theorem we get:

\subsection{Continued fractions for $\pi$}

The values $k=1, \, 3, \, 5$ give

\[
\frac{8}{4 -\pi}=3^2+ \cfrac{1^2 3^2}{(5^2-1^2)+ \cfrac{3^2 5^2}{(7^2-3^2)+ \cfrac{5^2 7^2}{(9^2-5^2)+ \cfrac{7^2 9^2}{(11^2-7^2)+ \dotsb}}}}\ ,
\]
\[
\frac{576}{75\pi -224}=7^2+ \cfrac{1^2 7^2}{(9^2-1^2)+ \cfrac{3^2 9^2}{(11^2-3^2)+ \cfrac{5^2 11^2}{(13^2-5^2)+ \cfrac{7^2 13^2}{(11^2-7^2)+ \dotsb}}}}
\]
and
\[
\frac{25600}{20992 -6615\pi}=11^2+ \cfrac{1^2 11^2}{(13^2-1^2)+ \cfrac{3^2 13^2}{(15^2-3^2)+ \cfrac{5^2 15^2}{(17^2-5^2)+ \cfrac{7^2 17^2}{(19^2-7^2)+ \dotsb}}}} .
\]

\subsection{Continued fractions for G}

The values $k=0, \, 2, \, 4$ give

\[
\frac{1}{G}=1^2+ \cfrac{1^4}{(3^2-1^2)+ \cfrac{3^4}{(5^2-3^2)+ \cfrac{5^4}{(7^2-5^2)+ \cfrac{7^4}{(9^2-7^2)+ \dotsb}}}}\ ,
\]
(this is the same one as (\ref{eq:bowman}))

\[
\frac{64}{19 -18G}=5^2+ \cfrac{1^2 5^2}{(7^2-1^2)+ \cfrac{3^2 7^2}{(9^2-3^2)+ \cfrac{5^2 9^2}{(11^2-5^2)+ \cfrac{7^2 11^2}{(13^2-7^2)+ \dotsb}}}}
\]
and
\[
\frac{49152}{22050G -19595}=9^2+ \cfrac{1^2 9^2}{(11^2-1^2)+ \cfrac{3^2 11^2}{(13^2-3^2)+ \cfrac{5^2 13^2}{(15^2-5^2)+ \cfrac{7^2 15^2}{(17^2-7^2)+ \dotsb}}}}.
\]

\section{Other continued fractions for $G$}
While looking for continued fractions for $G$ related to (\ref{eq:bowman}) in the same way that Lange's continued fraction for $\pi$ is related to Brouncker's, we found these:
\begin{equation}
\frac{2^5}{6G -1}=7+ 3\cdot \cfrac{1^4}{3(3^2-1^2)+ \cfrac{3^4}{3(5^2-3^2)+ \cfrac{5^4}{3(7^2-5^2)+ \cfrac{7^4}{3(9^2-7^2)+ \dotsb}}}} \label{eq:G}
\end{equation}
\[
\frac{2^{13}}{82G -19}=145+ 41\cdot \cfrac{1^4}{5(3^2-1^2)+ \cfrac{3^4}{5(5^2-3^2)+ \cfrac{5^4}{5(7^2-5^2)+ \cfrac{7^4}{5(9^2-7^2)+ \dotsb}}}}
\]
\[
\frac{2^{17}}{882G -\frac{713}{3}}=229+ 49\cdot \cfrac{1^4}{7(3^2-1^2)+ \cfrac{3^4}{7(5^2-3^2)+ \cfrac{5^4}{7(7^2-5^2)+ \cfrac{7^4}{7(9^2-7^2)+ \dotsb}}}}
\]
and it goes on like this. We'll prove (\ref{eq:G}). This continued fraction is related to the following series:
\[
\sum_{n=1}^\infty \frac{(-1)^{n-1}}{(4(n-1)^2+3)(4n^2+3)(2n-1)^2}
\]
as can be checked using Euler's method to transform a series into a continued fraction. Now the sum of this series can be calculated like this:
\begin{align*}
& \sum_{n=1}^\infty \frac{(-1)^{n-1}}{(4(n-1)^2+3)(4n^2+3)(2n-1)^2} \\
& = \sum_{n=1}^\infty {(-1)^{n-1}}\left(-\frac{1}{32 (4 (n-1)^2+3)}-\frac{1}{32 (4n^2+3)}+\frac{1}{16(2n-1)^2} \right)\\
& = -\frac{1}{32} \sum_{n=1}^\infty \frac{(-1)^{n-1}}{4 (n-1)^2+3} -\frac{1}{32} \sum_{n=1}^\infty \frac{(-1)^{n-1}}{4n^2+3} + \frac{1}{16}\sum_{n=1}^\infty \frac{(-1)^{n-1}}{(2n-1)^2}\\
& = -\frac{1}{3\cdot 32}+\frac{1}{16}G
\end{align*}
as a consequence of telescoping and the definition of $G$.\\
The series needed to prove the value of the other 2 continued fractions for $G$ are:
\[
\sum_{n=1}^\infty \frac{(-1)^{n-1}}{r(n-1) r(n) (2n-1)^2}
\]
with $r(n) = 16n^4+88n^2+41$ and $r(n)=64n^6+1168n^4+3628n^2+1323$ respectively.

\begin{center} 
{\textit{\large Appendix}}
\end{center}

\noindent\textit{The series $y_m$ ($m=1, 2, 3, \ldots$) defined by
\[
y_m = \sum_{n=1}^\infty \frac{(-1)^{n-1}}{(2n-1)^2(2n+1)^2(2n+3)^2 \cdots (2n+2m-3)^2}
\]
satisfy the recurrence relation:
\[
y_m=\frac{10m^2+8m+1}{2(2m+1)!!^2} -4m(m+1)^3 y_{m+2} .
\]}

\begin{proof}
By manipulating the terms of the series, we immediately get the result. We have:
{\small
\begin{align*}
y_{m+2} & = \sum_{n=1}^\infty \frac{(-1)^{n-1}}{(2n-1)^2 (2n+1)^2 \cdots (2n+2m-1)^2(2n+2m+1)^2}\\
& = \frac{1}{2m+2} \sum_{n=1}^\infty (-1)^{n-1} \left[\frac{1}{(2n-1)^2 (2n+1)^2 \cdots (2n+2m-1)^2(2n+2m+1)} \right.\\
& ~ \hspace{4cm} \left. -\frac{1}{(2n-1)(2n+1)^2 \cdots (2n+2m+1)^2}\right]\\
& =\frac{1}{(2m+2)^2}  \sum_{n=1}^\infty (-1)^{n-1} \left[\frac{1}{(2n-1)^2 (2n+1)^2\cdots (2n+2m-1)^2} \right.\\
& ~ \hspace{4cm}-\frac{1}{(2n-1)(2n+1)^2 \cdots (2n+2m-1)^2(2n+2m+1)}\\
& ~ \hspace{4cm}-\frac{1}{(2n-1)(2n+1)^2 \cdots (2n+2m-1)^2(2n+2m+1)}\\
& ~ \hspace{4cm}\left. +\frac{1}{(2n+1)^2(2n+3)^2 \cdots (2n+2m+1)^2}\right]
\end{align*}}\\
or
{\small
\begin{align*}
(2m+2)^2 y_{m+2} & =  \sum_{n=1}^\infty \frac{ (-1)^{n-1}}{(2n-1)^2 (2n+1)^2\cdots (2n+2m-1)^2} \\
& \hspace{1cm} - 2\sum_{n=1}^\infty \frac{ (-1)^{n-1}}{(2n-1)(2n+1)^2 \cdots (2n+2m-1)^2(2n+2m+1)} \\
& \hspace{1cm} -\sum_{n=2}^\infty \frac{ (-1)^{n-1}}{(2n-1)^2(2n+1)^2 \cdots (2n+2m-1)^2}
\end{align*}}\\
where we have shifted the summation index of the last series. Now the first and last sum cancel out, leaving one term, and the previous equation reduces to:
{\small
\begin{align*}
(2m+2)^2 y_{m+2}  =  \frac{1}{(2m+1)!!^2}-2 \sum_{n=1}^\infty \frac{(-1)^{n-1}}{(2n-1)(2n+1)^2 \cdots (2n+2m-1)^2(2n+2m+1)} .
\end{align*}}

We now repeat the two steps we used at the beginning of the proof to rewrite the remaining series. We get:
{\small
\begin{align*}
& \sum_{n=1}^\infty \frac{(-1)^{n-1}}{(2n-1)(2n+1)^2 \cdots (2n+2m-1)^2(2n+2m+1)} \\
& = \frac{1}{2m(2m+2)} \left[\frac{1}{ (2m-1)!!^2 (2m+1)}
- 2 \sum_{n=1}^\infty \frac{ (-1)^{n-1}}{(2n+1)^2 (2n+3)^2\cdots (2n+2m-1)^2} \right] \\
& = \frac{1}{2m(2m+2)} \left[\frac{1}{ (2m-1)!!^2 (2m+1)}
+ 2 y_m - \frac{2}{(2m-1)!!^2} \right]
\end{align*}}\\
where again we have shifted the summation index in the last step. Combining everything leads to the desired result.
\end{proof}

\noindent{\textbf{Note}:} E. Fabry (\textit{Th\'eorie des S\'eries \`a Termes Constants: Applications aux Calculs Num\'eriques,} A. Hermann \& fils, Paris, 1910, p.135) proved this recurrence using Kummer's transformation. Our proof is more straightforward.

\end{document}